\documentclass[11pt]{article}
\usepackage{diagrams}
\begin{document}
\title{\textbf{A note on smoothness and differential bases in positive characteristic}}
\author{Cristodor Ionescu\footnote{Work supported by the CNCSIS Grant GR 23/2008/812.}\\
Institute of Mathematics \textit{Simion Stoilow}\\ of the Romanian Academy,\\ P.O. Box 1-764, RO-014700, Bucarest, Romania,\\ 
e-mail: cristodor.ionescu@imar.ro}                           
\date{}                                                               
\newtheorem{theorem}{Theorem}[section]
\newtheorem{lemma}[theorem]{Lemma}
\newtheorem{proposition}[theorem]{Proposition}
\newtheorem{corollary}[theorem]{Corollary}
\newtheorem{definition}[theorem]{Definition}
\newtheorem{remark}[theorem]{Remark}
\newtheorem{example}[theorem]{ Example}
\newcommand{\dem}{\par\noindent{\it Proof: }}
\maketitle
\begin{abstract}
Let $u:A\to B$ be a morphism of noetherian local rings.  We  obtain  smoothness criteria for algebras with differential bases, in the  case of rings containing a field of characteristic $p>0.$ We also give smoothness criteria for reduced morphisms.
\textbf{Abbreviated title}: Smoothness and differential bases. 
\end{abstract}
\section{Introduction}
\par Let $p$ be a prime number. All over this paper, by a ring of characteristic $p$, we will mean a ring containing a field of characteristic $p.$ For our results we need some preliminary considerations. Let $A$ be a ring of  characteristic $p>0$. Denote by $A^{(p)}$ the $A$-algebra $A$ given by the Frobenius endomorphism of $A.$ For a ring morphism $u:A\longrightarrow B$ we have the following commutative diagram:

\begin{diagram}[size=2.5em]
      A             &   \rTo^{u}      & B                  &      &    \\
      \dTo<{F_A}    &                 & \dTo>{F_A\otimes 1_B}               &      & \\
      A^{(p)}       & \rTo     &  A^{(p)}\otimes_A B      &\rTo^{\omega_{B/A}} &B^{(p)}
      \end{diagram}
where $F_A$ is the Frobenius morphism of $A$ and 
$$\omega_{B/A}(a\otimes b)=u(a)\cdot b^p,\ \forall a\in A, b\in B.$$
\begin{remark}\label{0.0}
{\rm We shall freely use the definitions and notations from   \cite{EGA} and \cite{M} with respect to smoothness, regular  morphisms, p-basis etc.}
\end{remark}
\begin{definition}\label{1}
Let $u:A\rightarrow B$ be a morphism of noetherian rings. A system of elements $\{b_i\}_{i\in I}$ from $B$ is called a differential basis of $B$ over $A$, if the set $\{d_{B/A}(b_i)\}_{i\in I}$ is a basis of the module of differentials $\Omega_{B/A},$ where $d_{B/A}$ is the universal derivation of $B$ over $A.$
\end{definition}
\begin{remark}\label{2} {\rm 1) It is well-known that if $B$ has a $p$-basis  over $A,$ this is also a differential basis of $B$ over $A.$
\par\noindent 2) Obviously, if $B$ has a differential basis over $A$, then $\Omega_{B/A}$ is a free $B$-module.
\par\noindent 3)  Even when $A$ is a field, if  $\Omega_{B/A}$ is a free $B$-module,  $B$ does not have necessarily a differential basis over $A$. Indeed, let $A:={\bf F}_5,$ the prime field with 5 elements. Then such a simple ring as $B:={\bf F}_5[X,Y]/(X^2+Y^2-1),$ doesn't have a differential basis over $A$, while $\Omega_{B/A}$ is a free $B$-module \cite{T}.}
\end{remark}
We need the following important result of Ty\c c \cite[Theorem 1]{T}.
\begin{theorem}\label{22}
Let A be a noetherian ring and B a noetherian A-algebra containing A. Then any differential basis of B over A is a p-basis of B over A.
\end{theorem}
\dem See \cite[Prop. 58]{A} for a complete proof.

\section{A generalization of a result of Ty\c c}

In    \cite{T}, as a consequence of Theorem \ref{22}, A. Ty\c c proved the following smoothness criterion for a $k$-algebra, where $k$ is a field of characteristic $p>0.$  
\begin{theorem}\label{9}{\rm \cite[Theorem 2]{T}} 
Let $k$ be a field  of characteristic $p>0$ and $B$ a noetherian $k$-algebra. Suppose that:
\par\noindent a) B has a differential basis over k;
 \par\noindent b) $B\otimes_kk^{p^{-1}}$ is reduced.
 \par\noindent Then $B$ is smooth over $k$.
\end{theorem}
We give a generalization of Ty\c c's result, dropping the assumption that $k$ is a field. Our proof is slightly different from Ty\c c's proof.
 
\begin{theorem}\label{8}
Let $u:A\longrightarrow B$ be a  morphism of noetherian rings of characteristic $p>0.$ Suppose that:
\par\noindent a) $B$ has a $p$-basis over $A;$
\par\noindent b) $\omega_{B/A}$ is injective.
 \par\noindent Then $u$ is smooth.
\end{theorem}

\par\dem From a) it follows that $\Omega_{B/A}$ is a free $B$-module. Then it is enough to show that $u$ is regular. By \cite[Th\'eor\`eme 4]{Ra}, we have  to show that $\omega_{B/A}$ is flat. For this it is enough to find a $p$-basis of $B$ over $A^{(p)}\otimes_AB.$ Let $\{x_i\}_{i\in I}$ be a $p$-basis of $B$ over $A$. Since $\omega_{B/A}$ is injective we have that $A[B^p]\cong A^{(p)}\otimes_AB$ and  clearly 
 $\{x_i\}_{i\in I}$ is also a  $p$-basis  of $B$ over $ A^{(p)}\otimes_AB.$ 

\begin{corollary}\label{1000}
Let $u:A\longrightarrow B$ be an injective  morphism of noetherian rings of characteristic $p>0.$ Suppose that:
\par\noindent a) $B$ has a differential basis over $A;$
\par\noindent b) $A^{(p)}\otimes_A B$ is a reduced ring.
 \par\noindent Then $u$ is smooth. 
\end{corollary}

\par\dem Since $u$ is injective, by \ref{22} it follows that any differential basis of $B$ over $A$ is also a $p$-basis of $B$ over $A$. On the other hand, from the reducedness of $A^{(p)}\otimes_A B$ we get  that the Frobenius morphism of $A^{(p)}\otimes_A B,$ $F:A^{(p)}\otimes_A B\to (A^{(p)}\otimes_A B)^{(p)}$ is injective. But $F=(F_A\otimes_AB)^{(p)}\omega_{B/A},$ where
$$A^{(p)}\otimes_A B\stackrel{\omega_{B/A}}{\rightarrow} B^{(p)}\cong A\otimes_A B^{(p)}\stackrel{(F_A\otimes_AB)^{(p)}}{\rightarrow}(A^{(p)}\otimes_A B)^{(p)},$$
and $F(a\otimes b)=a^p\otimes b^p.$ It follows that $\omega_{B/A}$ is injective. Now we apply \ref{8}.

 \vspace{0.3cm}
 
 As a consequence of \ref{1000} we obtain the original result of Ty\c{c} (see \ref{9}).
 
 \vspace{0.3cm}

 \section{Smoothness and reduced morphisms}
In this section, mainly as consequences of the results in the previous section, we obtain some smoothness results for reduced morphisms. Let us first remind the definition of reduced morphisms.

\begin{definition}\label{5000}
 A morphism of noetherian rings $u:A\to B$ is called a reduced morphism if $u$ is flat and for every $p\in {\rm Spec}(A)$ and every field $L$, finite extension of $k(p)$, the ring $B\otimes_A L$ is a reduced ring. 
\end{definition}

\begin{theorem}\label{77} Let $u:A\longrightarrow B$ be an injective morphism of noetherian rings of characteristic p. Suppose that:
\par\noindent a) B has a differential basis over A;
\par\noindent b) $u$ is a reduced morphism.
\par\noindent Then u is smooth.
\end{theorem} 
\dem From \cite[Theorem 4]{D} it follows that $\omega_{B/A}$ is injective. Now we apply \ref{8} and \ref{22}.

\vspace{0.3cm}

\begin{corollary}\label{23.7} Let A be a  ring  of characteristic $p>0,$ with geometrically reduced formal fibers, $I$ an ideal of $A$ contained in the Jacobson radical and $B$ the completion of $A$ in the $I$-adic topology. If $B$ has a differential basis over A, then $B$ is smooth over A.
\end{corollary}

\vspace{0.3cm}
\begin{corollary}\label{23} Let A be a Nagata local ring of characteristic $p>0,$ $I$ an ideal of $A$  and $B$ the completion of $A$ in the $I$-adic topology. Assume that $B$ has a differential basis over A. Then $B$ is smooth over A. 
\par\noindent In particular, suppose that   $\widehat{A}$, the completion of $A$ in the topology of the maximal ideal, has a differential basis over $A$.  Then  $A$ is quasi-excellent.
\end{corollary}

\vspace{0.3cm}
\begin{corollary}\label{231} Let A be a  ring of characteristic $p>0.$ Assume that there exists $m\geq 1$ such that $A[[X_1,...,X_m]]$ has a differential basis over A.  Then $A[[X_1,...,X_n]]$ is smooth over A, for any $n\geq 1.$ In particular, $A$ has geometrically regular formal fibers.
\end{corollary}

\dem From \cite[Corollary 2.3]{Sh} we get that $A\rightarrow A[[X_1,...,X_m]]$ is a reduced morphism. The first assertion follows at once from Theorem \ref{77} and \cite[Theorem 2.2]{Ta1}. For the second one we use \cite[Theorem 3.3]{D1}.

\vspace{0.3cm}
\begin{remark}\label{obs}
{\rm Situations when $A[[X_1,...,X_n]]$ is smooth over $A$ are investigated in \cite{Ta1} and \cite{Ta}.}
\end{remark}
\begin{theorem}\label{577} Let $(A,m,k)$ be a noetherian local ring of characteristic p. The following are equivalent:
\par\noindent a) A has geometrically reduced formal fibers and  $[k:k^p]<\infty$;
\par\noindent b) A has geometrically regular formal fibers and  $[k:k^p]<\infty$;
\par\noindent c) There exists $m\geq 1$ such that $A[[X]]=A[[X_1,...,X_m]]$ has a $p$-basis over $A;$
\par\noindent d) There exists $m\geq 1$ such that $A[[X]]=A[[X_1,...,X_m]]$ has a differential basis over $A.$
\end{theorem} 
\dem It is well-known that a) and b) are equivalent. Also it follows from \ref{22} that c) and d) are equivalent. 
\par\noindent Let us prove that $b)\Rightarrow c)$. Let $B:=A[[X]]=A[[X_1,\dots,X_m]].$ Since $[k:k^p]<\infty$ we have that $A[B^p]=A^p[[X^p]][A]=A[[X^p]].$ Hence clearly $X$ is a $p-$basis of $A[[X]]$ over $A$.
\par\noindent In order to prove $d)\Rightarrow b)$ assume that $A[[X_1,...,X_m]]$ has a differential basis over $A$. From  Corollary \ref{231} we get that $A[[X_1,...,X_m]]$ is smooth over $A.$ Now we apply \cite[Theorem 2.2]{Ta1}.
  
\begin{remark}\label{1111} {\rm One can also prove $b)\Rightarrow d)$ directly as follows.
From b) it follows that $A$ is a finite $A^p$-algebra and by \cite[Theorem 2.2]{Ta1} we get that $B:=A[[X]]=A[[X_1,\dots,X_m]]$ is smooth over $A.$ We consider the $X-$adic topologies. Then $\Omega_{B/A}$ is a projective $B-$module and consequently  separated. Hence $\Omega_{B/A}\subseteq\widehat{\Omega_{B/A}}.$ On the other hand $\widehat{\Omega_{B/A}}$ is a free $B$-module, having as basis $dX\in\Omega_{B/A},$ where $d=d_{B/A}$ is the universal derivation. It follows that $\Omega_{B/A}$ is a free $B$-module having $dX$ as a basis.}
\end{remark}

\begin{remark}\label{obso}
{\rm In the situation of Theorem \ref{577}, if $A[[X_1,...,X_m]]$ has a differential basis (resp. $p$-basis) over $A$ for some $m\geq 1,$ it is clear that  $A[[X_1,...,X_n]]$ has a differential basis (resp. $p$-basis) over $A$ for all $n\geq 1.$}
\end{remark}
\vspace{0.3cm}
Finally we give an application of Theorem \ref{77} for integrable derivations \cite{MM}.
\begin{corollary}\label{23.1}
 Let $u:A\longrightarrow B$ be an injective morphism of noetherian rings of characteristic $p>0$. Suppose that:
\par\noindent a) B has a differential basis over A;
\par\noindent b) u is a reduced morphism.
\par\noindent Then ${\rm Der}_A(B)={\rm Ider}_A(B).$
\end{corollary}
\par\dem By Theorem \ref{77} $u$ is smooth. It follows   that $H^1(A,B,E)=0$, for every B-module E.   Now we apply \cite[Theorem 8]{MM}.

\vspace{0.3cm}

\begin{example}\label{777}
{\rm Let $k$ be a perfect field of characteristic 2 and let  
$$A:=(k[X,Y,Z]/(XY-Z^2))_{(X,Y,Z)}=k[x,y,z]_{(x,y,z)},\ xy=z^2.$$
Then $A$ is reduced, even normal, hence geometrically reduced over $k$. It is easy to see that  $\Omega_{A/k}$ is a free $A$-module with basis $\{xdx+ydy, dz\},$ where $d=d_{A/k}$ is the universal derivation. Since ${\rm Der}_k(A)\neq{\rm Ider}_k(A)$ \cite[Example 6]{MM}, it follows that $A$ doesn't have a differential basis over $k.$  }
\end{example}

\end{document}